\documentclass{article}

\usepackage{latexsym}
\usepackage{amssymb}
\usepackage{amsfonts}
\usepackage{amsmath}
\usepackage[all]{xy}
\usepackage{color}

\usepackage{mathrsfs}
\usepackage[utf8]{inputenc}
\usepackage[T1]{fontenc}
\usepackage{amscd}
\usepackage{amsxtra}

\newtheorem{thm}{\bf Theorem}[section]
\newtheorem{cor}[thm]{\bf Corollary}
\newtheorem{lem}[thm]{\bf Lemma}
\newtheorem{prop}[thm]{\bf Proposition}
\newtheorem{defn}[thm]{\bf Definition}

\newcommand{\field}[1]{\mathbb{#1}}

\newcommand{\N }{\field{N}}

\def\F{{\cal F}}
\def\A{{\cal A}}

\def\X{{\cal X}}
\def\Y{{\cal Y}}
\def\ZZ{{\cal Z}}
\newcommand{\HH}{{\cal H}}

\def\End{{\rm End}}
\def\Ext{{\rm Ext}}
\def\Hom{{\rm Hom}}
\def\Im{{\rm Im}}
\def\Mod{{\rm Mod}}
\def\proof{{\parindent0pt {\bf Proof.\ }}}

\def\coker{{\rm coker}}
\def\Ker{{\rm ker}}

\def\lfd{{\rm l_{\F}\!-\!dim}}
\def\rfd{{\rm r_{\F}\!-\!dim}}

\def\pd{{\rm proj.dim}}

\def\lxd{{\rm  l_{\X}\!-\!dim}}
\def\ryd{{\rm  r_{\Y}\!-\!dim}}

\def\elfd{{\rm e.l_{\F}\!-\!dim}}
\def\erfd{{\rm e.r_{\F}\!-\!dim}}
\def\elxd{{\rm e.l_{\X}\!-\!dim}}
\def\eryd{{\rm e.r_{\Y}\!-\!dim}}
\def\elfxd{{\rm e.l_{F(\X)}\!-\!dim}}
\def\elgxd{{\rm e.l_{G(\X)}\!-\!dim}}
\def\ergyd{{\rm e.r_{G(\Y)}\!-\!dim}}

\def\erfyd{{\rm e.r_{F(\Y)}\!-\!dim}}

\def\id{{\rm inj.dim}}

\def\rgldim{{\rm r.gl.dim}}
\def\lgldim{{\rm l.gl.dim}}
\def\lPCgldim{{\rm l.\mathscr{P}_C\!-\!gl.dim}}
\def\rICgldim{{\rm r.\mathscr{I}_C\!-\!gl.dim}}

\def\lXgldim{{\rm l.\X\!-\!gl.dim}}

\def\rYgldim{{\rm r.\Y\!-\!gl.dim}}

\newcommand{\cqfd}
{\hspace{1cm}
\rule{2mm}{2mm}%
\medbreak%
\par%
}

\begin{document}

\title{On proper and exact relative homological dimensions}

\author{Driss Bennis$^{1}$, J. R. Garc\'{\i}a Rozas$^2$, Lixin Mao$^3$ and Luis Oyonarte$^4$}

\date{}

\maketitle

\begin{abstract}
In Enochs' relative homological dimension theory occur the so called
(co)resolvent and  (co)proper dimensions which are defined using
proper and coproper resolutions constructed by  precovers and
preenvelopes, respectively. Recently, some authors have been
interested in relative homological dimensions defined by just exact
sequences.  In this paper, we contribute to the investigation of
these relative homological dimensions. We first study the relation
between these two kinds of relative homological dimensions and
establish  some ``transfer results"  under adjoint pairs. Then,
relative global dimensions are studied which lead to  nice
characterizations of some properties of particular cases of
self-orthogonal subcategories. At the end of the paper,  relative
derived functors are studied and generalizations of some known
results of  balance for relative homology are established.
\end{abstract}

\medskip
{\scriptsize 2010 Mathematics Subject Classification: Primary 18G25.
Secondary  16E30}


\section{Introduction}

\hskip .5cm Throughout this paper $R$ will be an associative (not
necessarily commutative) ring with identity, and all modules will
be, unless otherwise specified, unitary left $R$-modules.

We use $Proj(R)$ (resp., $Inj (R)$) to denote the class of all
projective (resp., injective) $R$-modules. The category of all left $R$-modules will be denoted by $R$-Mod.
For an $R$-module $C$, we use $Add_R(C)$   to denote the class of all $R$-modules which are
isomorphic to direct summands of direct sums  of copies of $C$, and $Prod_R(C)$ will denote the class of all $R$-modules
which are isomorphic to direct summands of direct products of copies of $C$. Also $\A$ will be an abelian category. By
a subcategory of a given category $\A$ we will always mean a full
subcategory closed under isomorphisms and having the zero object of
$\A$.
 \medskip

Let us start with some basic definitions.

Given a subcategory   $\F$ of  $\A$,  an $\F$-precover of an object
$M\in \A$ is a morphism $\varphi:F\to M$ ($F\in \F$) such that
$\Hom(F',\varphi)$ is epic for every $F'\in\F$. An $\F$-precover
$\varphi$ is said to be special if it is epic and $\Ext^1(F',\ker
\varphi)=0$ for every $F'\in \F$. The subcategory $\F$  of  $\A$ is
said to be precovering if every object  of  $\A$ has an
$\F$-precover. Dually, an $\F$-preenvelope of an object  $M\in \A$
is a morphism $\varphi: M \to F$ ($F\in \F$) such that $\Hom(\varphi
,F')$ is epic for every $F'\in\F$. An $\F$-preenvelope $\varphi$ is
said to be special if it is monic and $\Ext^1(\coker \varphi,F')=0$
for every $F'\in \F$. The subcategory $\F$ is said to be
preenveloping if every object has an $\F$-preenvelope.

Given a subcategory   $\F$ of  $\A$, a proper left $\F$-resolution
of an object $M$ of $\A$ is a complex (so not necessarily exact) of
objects in $\F$: $$\textbf{X}=\cdots  \rightarrow X_1 \rightarrow
X_0 \rightarrow M \rightarrow 0$$ such that $\Hom(F,\mathbf{X})$ is
an exact complex for every $F\in \F$. Thus we see that $M$ having a
proper left $\F$-resolution is equivalent to $M$ having an
$\F$-precover whose kernel has an $\F$-precover and so on. Coproper
right $\F$-resolutions may be defined dually.

An object $M$ of $\A$ is said to have proper left $\F$-dimension  less than or equal to an integer $n\geq 0$, $\lfd(M)\leq n$, if $M$ admits a proper left $\F$-resolution  of the form
$$0\to A_n\to\cdots\to A_1\to A_0\to M\to 0.$$  If $n$ is the least nonnegative integer for which such a resolution exists then we set $\lfd(M)=n$, and if there is no
such $n$ then we set $\lfd (M)=\infty$. Dually, the right
$\F$-dimension, $\rfd(M)$, is defined when $\F$ is preenveloping.
The reader is invited to see Enochs and Jenda's book \cite{EJb} for
more details.

It has been shown that (co)proper resolutions and dimensions are
suitable in the context of relative homological algebra in order to
establish a theory analogous to the classical one. So relative
derived functors were defined and various and interesting results
using the relative proper and coproper dimensions were proved.
However, it does not seem that this relative homological dimension
is suitable for measuring how far away a  module is in the given
class.  That is the reason why some authors prefer to work with the
following dimensions which are defined using  exact sequences.

\begin{defn}\label{def-exact-dim} Given a subcategory   $\F$ of  $\A$.
An object  $M$ of $\A$ is said to have exact left $\F$-dimension
less than or equal to $n$, $\elfd(M)\leq n$, if there is an exact
sequence
$$0\to A_n\to\cdots\to A_1\to A_0\to M\to 0$$ with $A_i\in \F$ for every $ i\in \{0,...,n\}$. If $n$ is the least nonnegative integer for which such a
sequence exists then $\elfd(M)=n$, and if there is no such $n$, then
$\elfd(M)=\infty$.

Dually the exact right $\F$-dimension, $\erfd(M)$,  is defined.
\end{defn}

But, using only exact sequences would not lead to interesting
results. That is the reason why some authors assume more conditions
on the class that help to build a new relative homological dimension
theory as given in the two papers \cite{Hu} and \cite{Zhu14}.

The discussion above is the reason behind the following declaration
of Takahashi and White: ``Because proper left PC-resolutions need
not be exact, it is not immediately clear how best to define the
left  PC-dimension  of a module. For instance, should one consider
arbitrary left  PC-resolutions or only exact ones?". After that they
proved that indeed  they are  the same for the classes of the
so-called ${\mathscr I}_C$-injective and  ${\mathscr
P}_C$-projective modules, where $C$ is a semidualizing module.
However, the relation between the two kinds of dimensions is still
an open question.

Our aim in this paper is  twofold, to study the relation between the
two kinds of dimensions and also to continue the investigation by
involving the relative derived functors.

We have structured the paper in the following way:

In Section 2, we investigate the relation between the two kinds of
dimensions. For this reason we introduce some new kinds of
subcategories on which these dimensions are closely related. Namely,
a subcategory $\X$ is said to be EP, if every object with finite
exact dimension relative to $\X$  has also finite proper
$\X$-dimension.  We introduce also PE subcategories when we are
interested in the converse implication. Among other results, we give
two important  situations of EP subcategories   (Theorems
\ref{thm-wEP} and \ref{thm-Cover-EP}). We end Section 2 with a
result that studies PE subcategories (see Theorem \ref{thm-eld-ld}).

Section 3 is devoted to   some ``transfer results"  under adjoint pairs
(see Theorems
\ref{thm-compa-1} and \ref{thm-compa2}).

In Section 4, we investigate the  global dimensions relative to
subcategories. We start with a general study (Proposition
\ref{prop-gldim}) and then we investigate some interesting
particular cases (see Corollary \ref{cor-gldim-2} and Propositions
\ref{prop-gldim-sigma} and \ref{prop-gldim-Gen}).

The last section of the paper is interested in the relative
(co)homology groups. The first main result  of this section
characterizes relative dimensions using relative cohomology groups
(Theorem \ref{thm-rel-func}). The second aim is the   balance
results for relative (co)homology groups (see Theorems
\ref{thm-compa-3} and \ref{thm-compa-5}).

\section{EP and PE subcategories}

\hskip .5cm  Our aim in this section is to investigate the relation  between
exact dimensions and proper (coproper) dimensions. To this end we
introduce the following notions:

\begin{defn}\label{def-ep6pe}\begin{enumerate}
\item A subcategory $\F$ of  $\A$ is called
 EP   if, for any object  $M$ of $\A$, if $\elfd(M)\leq n$ (for some integer $n$), then  $\lfd(M)\leq n$.
 \item A subcategory $\F$ of  $\A$ is called
 PE   if, for any object  $M$ of $\A$, if $\lfd(M)\leq n$ (for some integer $n$), then  $\elfd(M)\leq n$.
\end{enumerate}
\end{defn}

Of course, one can define as above EC and   CE subcategories using
(coproper, exact) right dimensions. All the results given in this
section have dual ones using EC and CE terminologies. Hence we do
not need to state these dual results except  for Theorem
\ref{thm-eld-ld} (see Theorem \ref{thm-erd-rd}).

In terms of the definition above, our aim in this section is   to
study EP and PE subcategories.

Let us first notice that by ``$\F$ is EP"  it should not be
understood  that, for any object  $M$ of $\A$, if $\elfd(M)\leq n$,
then there is an exact proper left resolution of $M$ of length $n$.
In fact, this would be another interesting question. When $\F$ is
 closed under extensions and direct summands we have a partial positive answer for an
object  $M$ of $\A$ with $\elfd(M)\leq 1$. Indeed, in this case, $M$
admits both an exact sequence $0 \to K''  \to  K' \to M\to 0$, with
$K'$ and $K''$ in $\F$, and an $\F$-precover   $P \to M\to 0$. Let $
K=\ker (   P \to M)$ and consider the following diagram:
$$\xymatrix{ 0\ar[r]& K''\ar[r]   \ar[d] &  K' \ar[d] \ar[r] & M \ar@{=}[d]  \ar[r] & 0\\
0\ar[r]& K \ar[r]&  P    \ar[r]& M\ar[r]   & 0 }  $$
which induces a mapping cone diagram
$$ \xymatrix{ 0\ar[r]& K''   \ar[r]    &  K' \oplus K    \ar[r] & P   \ar[r] & 0.}  $$
Therefore, $K$ is also in $\F$ since $\F$ is closed under extensions
and direct summands.

The following result gives  a partial positive answer for  a
precovering class.

\begin{thm}\label{thm-Cover-EP}
Assume that $\X$ is a precovering class  and closed under extensions
and direct summands. If any object of $\A$ has an epic
$\X$-precover, then, for any object  $M$ of $\A$ with $\elxd(M)\leq
n$ (for some integer $n$), there exists  an exact proper left
$\X$-resolution  of the form $$0\to A_n\to\cdots\to A_1\to A_0\to
M\to 0.$$ In particular, $\X$ is an EP subcategory of $\A$.
\end{thm}
\proof  Let $M$ be an object of $\A$ with $\elxd(M)\leq n$ (for some
integer $n$). For $n=0$ there is nothing to prove. So we can suppose
that $n>0$. Since  $\elxd(M)\leq n$, there is an exact sequence
$$0\rightarrow K_1  \rightarrow X_{0}\rightarrow  M \rightarrow 0$$
where $X_0\in \X$ and $\elxd(K_1)\leq n-1$. We prove, by induction
on $n\geq 1$,  that, for every $\X$-precover   $ Y_0 \to M\to 0$, $
K=\ker (   Y_0 \to M)$ has an exact proper left $\X$-resolution of
length $n-1$. This gives the desired result. For the case $n=1$, see
the discussion before Theorem \ref{thm-Cover-EP}. Now, for $n>1$,
consider  an $\X$-precover $A_0\to K$ of $K$ and the following
diagram which holds by Horseshoe lemma:
$$\xymatrix{
  &  0 \ar[d] & 0 \ar[d]  &  0 \ar[d]  &  \\
 0\ar[r]& K'\ar[r] \ar[d]  & L \ar[d] \ar[r] & K_1 \ar[d]  \ar[r] & 0\\
 0\ar[r]& A_0\ar[r] \ar[d]  & A_0 \oplus X_0 \ar[d] \ar[r] & X_0 \ar[d]  \ar[r] & 0\\
0\ar[r]& K \ar[d] \ar[r]&   Y_0   \ar[d] \ar[r]& M\ar[r]  \ar[d]& 0\\
  &  0&0   &0  &   }  $$
where  $ K'=\ker (  A_0\to K)$. It is clear that the middle vertical
sequence is $\Hom(\X,-)$-exact which, by applying $\Hom(X_0,-)$,
implies that it splits. So, $L$ is also in $\X$  (using the middle
horizontal sequence). On the other hand, the top horizontal sequence
is also  $\Hom(\X,-)$-exact which shows that $L\to K_1$ is an epic
$\X$-precover of $K_1$. Since $\elxd(K_1)\leq n-1$ and by induction,
$K'$  has an exact proper left $\X$-resolution of length $n-2$. So
using the left vertical sequence we get the desired result.\cqfd

\begin{cor}\label{cor-flat-EP} The class of flat modules is EP.
\end{cor}

The following gives another example  of EP subcategories.

We will use the following notations. Given a subcategory $\F$ of
$\A$, the subcategory of all object $N$ such that $\Ext^{\geq
1}_R(F,N)=0 $ for every $ F\in \F$ will be denoted by $\F^{\perp}$
(similarly, $^{\perp}\F=\{ N:\ \Ext^{\geq 1}_R(N,F)=0\ \forall F\in
\F\}$).

Recall that  a subcategory  $\HH$ of $\F$ is said to be an
Ext-injective cogenerator for $\F$ if $\HH \subseteq \F^{\perp}$
and, for each object $M\in\F$, there exists an exact sequence $0 \to
M \to  H \to M'\to 0$ such that $H\in \HH$  and $M'\in \F$ (see
\cite{Zhu14}).   In \cite{Zhu14}, various examples of classes which
admit Ext-injective cogenerators were given. For example, the class
of Gorenstein projective modules admits the class of projective
modules as an  Ext-injective cogenerator subcategory (See also \cite[Remark 2.7]{BGO1} for a more general context).

\begin{thm}\label{thm-wEP}
Any subcategory $\F$ of $\A$ which is closed under extensions and
admits an  Ext-injective cogenerator  $\HH$ is an EP subcategory.
\end{thm}
\proof  Let $M$ be an object of $\A$ with $\elfd(M)\leq n$ (for some
integer $n$). We claim that   $\lfd(M)\leq n$. Consider an exact
sequence
 $$0\rightarrow X_{n} \rightarrow  X_{n-1}\rightarrow \cdots \rightarrow X_1
\rightarrow X_0 \rightarrow M \rightarrow 0$$ where $X_i\in \F$. We
split it into short exact sequences $0\rightarrow K_{i} \rightarrow
X_{i-1}\rightarrow  K_{i-1}\rightarrow 0$, where $ K_{n}= X_{n}$, $
K_{i}=\Ker ( X_{i-1}\rightarrow  X_{i-2})$ for $i=1,...,n-1$ with $
X_{-1}= M$. Since $\HH$ is an  Ext-injective cogenerator, there
exists an exact sequence $0 \to  K_{i}  \to  H_i \to L_i\to 0$ such
that $H_i\in \HH$  and $L_i\in \F$. Then, we get the following
pushout diagram:
$$\xymatrix{ & 0 \ar[d] & 0 \ar[d]  &   &  \\
0\ar[r]&  K_{i}  \ar[d] \ar[r] & X_{i-1} \ar@{-->}[d] \ar[r] & K_{i-1} \ar@{=}[d]
\ar[r] &
0\\
0\ar[r]& H_i\ar[d]  \ar@{-->}[r] & D_i\ar[d] \ar[r] &  K_{i-1} \ar[r] & 0\\
 &  L_i \ar[d] \ar@{=}[r] &L_i \ar[d]  &   &  \\
     &   0 & 0   &   & }$$
If $n=1$, the middle horizontal sequence shows, as desired, that also
$\lfd(M)\leq 1$ since $\HH \subseteq \F^{\perp}$. If $n>1$, we continue our process. Since $\F$  is closed under extensions, $ D_{i}$ is in $\F$, and so we can consider
an exact sequence $0 \to  D_{i}  \to  H_{i-1} \to N_{i-1}\to 0$ such that $H_{i-1}\in \HH$  and $N_{i-1}\in \F$ and another   pushout diagram:
$$\xymatrix{
   &   &   0 \ar[d] & 0 \ar[d]  &    \\
0\ar[r]&  H_{i}  \ar@{=}[d] \ar[r] & D_{i} \ar[d] \ar[r] & K_{i-1} \ar@{-->}[d]
\ar[r] &
0\\
0\ar[r]& H_i\ar[r] & H_{i-1}\ar[d] \ar@{-->}[r]   &  E_{i-1} \ar[d]\ar[r] & 0\\
    &  &  N_{i-1} \ar[d] \ar@{=}[r] &N_{i-1} \ar[d] &  \\
   &  &   0 & 0   &    }$$
Then, using the  middle horizontal sequence, $ E_{i-1}$ is in $\F^{\perp}$. Consider now the  the following
 pushout diagram:
$$\xymatrix{
     &   0 \ar[d] & 0 \ar[d]  &   &  \\
0\ar[r]&  K_{i-1} \ar[d] \ar[r] & X_{i-2} \ar@{-->}[d] \ar[r] & K_{i-2} \ar@{=}[d]
\ar[r] &
0\\
0\ar[r]& E_{i-1}\ar[d]  \ar@{-->}[r] & D_{i-1}\ar[d] \ar[r] &  K_{i-2} \ar[r] & 0\\
 &   N_{i-1}  \ar[d] \ar@{=}[r] & N_{i-1} \ar[d]  &   &  \\
     &   0 & 0   &   & }$$
Therefore, the following exact sequence
$$\xymatrix{
 0\ar[r]& H_{i} \ar[r]& H_{i-1} \ar[r] & D_{i-1}  \ar[r] &  K_{i-2} \ar[r] & 0 }$$
and the process above lead to  the desired sequence.
 \cqfd

One can also be interested in the following  ``strong" condition on
EP subcategories: a subcategory $\F$ of  $\A$ is said to be strongly
EP, if any   exact sequence of this form $$0\to A_n\to\cdots\to
A_1\to A_0\to M\to 0$$ with each $A_i$ in  $\F$, is
$\Hom(\F,-)$-exact.

Clearly every  strongly EP subcategory is EP. We will give an example which shows that the converse is not true in general.
\begin{prop}\label{prop-elxd} A  subcategory $\X$ of $\A$, which is closed under extensions, is  strongly EP if and only if it is self-orthogonal (that is $\X\subseteq  \X^\perp$).
In this case, every object of finite   exact left $\X$-dimension belongs to the class $\X^\perp$.
\end{prop}
\proof  $\Rightarrow $.
Consider two objects $A$ and $C$ in $\X$. We claim that every short exact  sequence
$$   \xymatrix{
 0\ar[r]& A \ar[r] & B \ar[r]& C \ar[r] & 0}  $$
splits. Indeed,    $B$ is also in $\X$ (by hypothesis) and since $\X$ is    EP, the above sequence is $\Hom(\X,-)$-exact, in particular,
$$   \xymatrix{
 0\ar[r]& \Hom(C,A) \ar[r] & \Hom(C,B)  \ar[r]& \Hom(C,C) \ar[r] & 0}  $$
is exact, as desired.\\
 $\Leftarrow $. Consider an exact sequence of the form $$0\to
A_n\to\cdots\to A_1\to A_0\to M\to 0$$ with each $A_i$ in  $\X$. The case $n=0$ is trivial. Then assume that $n\geq 1$. Consider
the short exact  sequences $0\rightarrow K_i \rightarrow  A_{i-1} \rightarrow K_{i-1} \rightarrow 0$, where
 $K_{i}=\ker(A_{i-1}   \rightarrow  A_{i-2}) $. Using the fact that $A_i\in \X^{\perp}$, an induction process starting
 from $K_n=A_n \in \X^{\perp}$ shows that $K_{i}\in \X^{\perp}$ for every $  i$. Therefore,  the sequence is  indeed a proper left $\X$-resolution of $M$.

Finally, the last statement is obtained by applying the last argument to the sequence   $0\rightarrow K_1 \rightarrow  A_{0} \rightarrow M \rightarrow 0$.\cqfd
\medskip

Now it is clear that not every EP subcategory is a strongly EP
subcategory. For this consider the subcategory  $\F(R)$ of flat
modules  over a  ring $R$, then $\F(R)$ is  EP (because it is a
covering class) but it is not   strongly EP (because it is closed
under extensions but not necessarily self-orthogonal).

The exact dimensions relative to self-orthogonal (strongly EP)
subcategories behave similarly as the classical dimensions. Namely
we have the following result  which will be used in the proof of the next
theorem.

\begin{prop}\label{prop-elxd2} Assume $\X$ to be a self-orthogonal subcategory of $\A$ which is closed under extensions and direct summands.
Then, for   every   object $M$ of $\A$ with $\elxd( M)\leq n$ ($n\in
\N$), every   proper left $\X$-resolution $$\cdots \rightarrow Y_1
\stackrel{f_1}{\longrightarrow} Y_0 \stackrel{f_0}{\longrightarrow}
M\stackrel{f_{-1}} \longrightarrow 0$$ of $M$ is  exact and
$\ker(f_{n-1})\in \X$. In particular,  $M$ has a special
$\X$-precover.
\end{prop}
\proof Consider a left proper $\X$-resolution $$\cdots  \rightarrow Y_1 \stackrel{f_1}{\longrightarrow} Y_0 \stackrel{f_0}{\longrightarrow}
M\stackrel{f_{-1}} \longrightarrow 0$$ of $M$. We may assume that $n\geq 1$. So  $f_0$ is epic and $K=\ker(f_0)\in \X^\perp$. We  prove that
$\elxd( K)=n-1$ and this gives the desired result.
Since $\elxd( M)=n$, there exists an exact proper left $\X$-resolution of $M$ of the form $$0\rightarrow X_n \rightarrow \cdots  \rightarrow X_1 \rightarrow X_0 \rightarrow M \rightarrow 0.$$
  Let  $K_0=\Ker(  X_0 \rightarrow M )$.
By the relative version of Schanuel's lemma   \cite[Lemma 8.6.3]{EJb}, $ K \oplus X_0 \cong K_0 \oplus Y_0 $.
  Then, if $n=1$ ($K_0= X_1$), we get $K\in \X$, as desired. Now suppose that $n>1$ and consider the short exact sequences
$$0\rightarrow K  \rightarrow K_0\oplus  Y_0 \rightarrow X_0 \rightarrow 0$$
and $$0\rightarrow K_1 \rightarrow   X_1 \rightarrow K_0 \rightarrow 0,$$ where $K_1=\ker(  X_1 \rightarrow  X_0 )$. So in the following pullback diagram
$$   \xymatrix{
     &  0 \ar[d]  & 0 \ar[d]  &   &  \\
 &  K_1\ar[d] \ar@{=}[r]&  K_1\ar[d]  &   &  \\
 0\ar[r]& D\ar@{-->}[d] \ar@{-->}[r] & X_1 \oplus  Y_0\ar[d] \ar[r] &  X_0 \ar@{=}[d]  \ar[r] & 0\\
0\ar[r]& K\ar[r] \ar[d]&  K_0\oplus  Y_0  \ar[d] \ar[r]& X_0\ar[r] & 0\\
 & 0 &0  &   & }  $$
  the left vertical sequence gives $D\in  \X^\perp$. Then, the
middle horizontal sequence splits and so $D \in  \X$. Hence, the
left vertical sequence shows that $K$ admits an exact proper left
$\X$-resolution of the length $n-1$. Then, $\elxd(K_1)\leq n-1$. It
is clear that  $\elxd(K_1)< n-1$, which contradicts the fact that
$\elxd( M)=n$. Therefore,  $\elxd(K_1)= n-1$. \cqfd

Now we are interested in PE subcategories. The following   result characterizes when resolvent dimensions relative to self-orthogonal subcategories
are exact dimensions. For this, we need the following notions which extend those of  \cite[Definition 2.9]{BGO2}.

\begin{defn}
\begin{enumerate}
\item A subcategory $\X$ of $\A$ is said to be $\Hom$-faithful if for every object $N$ of $\A$: $\Hom(M,N)=0$
 for every $M \in \X$ implies that $N=0$.\item A subcategory $\Y$ of $\A$   is said to be $\Hom$-cofaithful  if for every object $M$ of $\A$: $\Hom (M,N)=0$ for every $N \in \Y$ implies that $M=0$.
\end{enumerate}
\end{defn}

The following result generalizes  \cite[Theorem 2.11]{BGO2} using  similar  arguments.  However, we give here a complete proof for reader's convenience.

\begin{thm}\label{thm-eld-ld} Assume $\X$ to be a self-orthogonal subcategory of $\A$ which is closed under extensions and direct summands.
The following assertions are equivalent.
\begin{enumerate}
\item $\X$ is $\Hom$-faithful.
\item  If $\phi: I\rightarrow M$ is an $\X$-precover with $K=\ker(\phi) \in \X^{\perp}$, then $\phi$ is epic and $M\in  \X^{\perp}$.
\item  If   $\cdots  \rightarrow Y_1 \stackrel{f_1}{\longrightarrow} Y_0 \stackrel{f_0}{\longrightarrow} M \longrightarrow 0$ is a proper left $\X$-resolution of an
object $M$ of $\A$ with $\ker(f_{n})\in
\X^{\perp}$ for some positive integer $n$, then the sequence $0\rightarrow \ker (f_n) \rightarrow \cdots  \rightarrow Y_1 \rightarrow Y_0 \rightarrow M \rightarrow 0$ is exact.
\item   Every  monic $\X$-precover $\phi: I\rightarrow M$  of an object $M$ of $\A$ is an isomorphism.
\item Every object $M$  of $\A$ of finite $\lxd (M)$  has a special $\X$-precover.
\item   $\X$ is a PE subcategory   of $\A$.
\item  If  $\lxd (M) <\infty$ for an object $M$ of $\A$, then  every left $\X$-resolution of $M$ is exact.
\end{enumerate}
\end{thm}
\proof $1.\Rightarrow 2$. Let $X \in \X$, since $\phi: I\rightarrow
M$ is an $\X$-precover, the sequence $$0\rightarrow \Hom(X,K)
\rightarrow \Hom(X,I) \rightarrow \Hom(X,M) \rightarrow 0$$ is
exact, and since $K\in \X^{\perp}$, we get the exact sequences
$$0\rightarrow \Hom(X,K) \rightarrow \Hom(X,I) \rightarrow \Hom(X,\Im (\phi)) \rightarrow 0$$ and
$$0=\Ext^{i}(X,I) \rightarrow \Ext^{i}(X,\Im (\phi))\rightarrow \Ext^{i+1}(X,K)=0\ \forall i\geq 1.$$ Thus,
$\Ext^{\geq 1}_R (X,\Im(\phi))= 0$ (so $\Im(\phi)\in \X^{\perp}$) and we have the following  commutative diagram:
$$\xymatrix{      &    &    & 0 \ar[d]  &  \\
     0\ar[r]  & \Hom(X,K) \ar@{=}[d] \ar[r]   &  \Hom(X,I) \ar@{=}[d] \ar[r]  & \Hom(X,\Im (\phi)) \ar[d] \ar[r]   & 0 \\
     0\ar[r]  & \Hom(X,K)   \ar[r]   &   \Hom(X,I)  \ar[r]  &  \Hom(X,M) \ar[d] \ar[r]   & 0 \\
       &    &   &  \Hom(X,M/\Im(\phi)) \ar[d]  &   \\
   &    &   &  0 &
}$$
This implies that $\Hom(X,M/\Im (\phi))=0$ and then that $M/\Im(\phi)=0$ since $\X$ is $\Hom$-faithful.

$2.\Rightarrow 1$. If $\Hom(X,M)=0$ for every $X\in \X$, then $0\to
M$ is a special $\X$-precover and so it is epic by hypothesis, hence
$M=0$.

$2.\Rightarrow 3$. This is proved by decomposing the sequence $$0\rightarrow \ker(f_n) \rightarrow \cdots  \rightarrow Y_1 \rightarrow Y_0 \rightarrow M \rightarrow 0$$
into short sequences and then applying
recursively the assertion 2.

$3.\Rightarrow 2$. The assertion 2 is a particular case of the assertion 3.

$2. \Rightarrow 4$. Clear.

$4.\Rightarrow 1$.   If $\Hom(X,M)=0$ for every $X\in \X$, then the
trivial homomorphism  $0\rightarrow M$ is a monic  $\X$-precover of
$M$. So it is an isomorphism which implies that $M=0$.

$3.\Rightarrow 5$. Use Proposition \ref{prop-elxd2}.

$5.\Rightarrow 4$. Consider  an object $M$ of $\A$ which admits a
monic $\X$-precover $\phi: I\rightarrow M$. Then, $\lxd (M)=0$. So
$M$ has a special $\X$-precover, in particular epic  $\X$-precover.
Then, $\phi$ is also epic, as desired.

$3.\Rightarrow 6$. Obvious.

The implications $6.\Rightarrow 5$, $7.\Rightarrow 5$ and $6.\Rightarrow 7$ are simple consequences of Proposition \ref{prop-elxd2}.\cqfd

Note that, in \cite[Definition 4.5]{H}, Holm used the term
``separating subcategory"  which we call here  ``$\Hom$-faithful
subcategory". In \cite[Lemma  4.6]{H}, he proved that, for any
precovering class $\F$, if    every monic $\F$-precover is an
isomorphism then $\F$ is separating. Here, we have proved that the
converse holds for self-orthogonal subcategories.

Dually we get the following result.

\begin{thm}\label{thm-erd-rd}  Assume $\Y$ to be a self-orthogonal subcategory of $\A$ which is closed under extensions and direct summands.
The following statements are equivalent.
\begin{enumerate}
\item $\Y$ is $\Hom$-cofaithful.
\item If $\phi:M\to Y$ is a $\Y$-preenvelope with $\coker (\phi)\in {^{\perp}\Y}$, then $\phi$ is injective and $M\in {^{\perp}\Y}$.
\item Given any coproper right $\Y$-resolution $0\to M\stackrel{f_0}{\to} Y_0\stackrel{f_1}{\to} Y_1\to \cdots $ of an object $M$ of $\A$, if $\Im (f_n)\in {^{\perp}\Y}$
for some positive integer $n$, then
the sequence $0\to M\to Y_0\to Y_1\to \cdots \to Y_{n-1}\to \Im (f_n)\to 0$ is exact.
\item Every epic $\Y$-preenvelope of any object of $\A$  is an isomorphism.
\item Every  object of $\A$  of finite right $\Y$-injective dimension has a special $\Y$-preenvelope.
\item $\Y$  is a CE subcategory   of $\A$.
\item  If   $\ryd(M)<\infty$ for an object $M$ of $\A$, then every right $\Y$-resolution of $M$ is exact.
\end{enumerate}
\end{thm}

As a functorial description of the modules of finite $\X$-projective
dimension, we have the following result.

\begin{cor}\label{cor-fun-ld} Under the hypotheses of Theorem \ref{thm-eld-ld},  if $\X$ is   $\Hom$-faithful, then, for every
$M\in \A$ of finite  left $\X$-dimension and every   nonnegative integer $n$,  the following assertions are equivalent.
\begin{enumerate}
\item  $\lxd (M)\leq n$.
\item  $\Ext^{>n} (M,X)=0$ for every $X\in \X$.
\item  $\Ext^{>n} (M,X)=0$ for every $X$ with finite left $\X$-dimension.
\end{enumerate}
\end{cor}
\proof   Follows from standard arguments.\cqfd

\medskip

Of course, the last result has its dual.

\begin{cor}\label{pro-fun-rd} Under the hypotheses of Theorem \ref{thm-erd-rd},  if $\Y$ is   $\Hom$-cofaithful, then, for every
$M\in \A$ of finite  right  $\Y$-dimension and every   nonnegative integer $n$,  the following assertions are equivalent.
\begin{enumerate}
\item  $\ryd(M)\leq n$.
\item  $\Ext^{>n} (Y,M)=0$ for every $Y\in \Y$.
\item  $\Ext^{>n} (Y,M)=0$ for every $Y $ with finite right $\Y$-dimension.
\end{enumerate}
\end{cor}

\section{Transfer results under adjoint pairs}

\hskip .5cm  We give in this section   some transfer results under adjoint pairs.
In particular, they generalize \cite[Theorem 2.11]{TW}.

We denote by $(F,G):  \mathcal{A}\rightarrow \mathcal{B}$  an
adjoint pair.  The natural transformations $\varepsilon: FG
\rightarrow 1_{\mathcal{B}}$ and $\eta: 1_{\mathcal{A}}\rightarrow
GF$ will mean the counit and the unit of the adjunction $(F,G)$.

The following lemma is a slight generalization of  \cite[Proposition 2.5]{GT1}.

\begin{lem}\label{lem-compa-1} Let $(F,G):  \mathcal{A}\rightarrow \mathcal{B}$ be an
adjoint pair with $\mathcal{Y}\subseteq\mathcal{A}$ and $\mathcal{X}\subseteq\mathcal{B}$. Then, $F(f): F(U)\rightarrow F(V)$ is an $F(\mathcal{Y})$-preenvelope
of $F(U)$ if $f: U\rightarrow V$ is a
$\mathcal{Y}$-preenvelope of $U$ in $\mathcal{A}$.
\end{lem}
 \proof Let $g:
F(U)\rightarrow F(W)$ be a morphism with $W\in \mathcal{Y}$, then we get the following commutative diagram:
$$\xymatrix{FGF(U)\ar[d]_{FG(g)} \ar[r]^{\varepsilon_{F(U)}}&F(U)\ar[d]_{g}\\
FGF(W)\ar[r]^{\varepsilon_{F(W)}}&F(W).}$$ Note that $GF(W)\in\mathcal{Y}$ and so there exists $h: V\rightarrow GF(W)$ such that $hf=G(g)\eta_{U}$. Thus we have
$$(\varepsilon_{F(W)}F(h))F(f)=\varepsilon_{F(W)}F(hf)=\varepsilon_{F(W)}F(G(g)\eta_{U})$$ $$=\varepsilon_{F(W)}FG(g)F(\eta_{U})=g\varepsilon_{F(U)}F(\eta_{U})=g.$$
Hence $F(f): F(U)\rightarrow F(V)$ is an $F(\mathcal{Y})$-preenvelope of $F(U)$.
\cqfd

\begin{thm}\label{thm-compa-1} Let $(F,G):  \mathcal{A}\rightarrow \mathcal{B}$ be an
adjoint pair with $\mathcal{Y}\subseteq\mathcal{A}$ and
$\mathcal{X}\subseteq\mathcal{B}$.
\begin{enumerate}\item If $FG(\mathcal{X})\subseteq\mathcal{X}$, and $\mathcal{D}\supseteq\mathcal{X}$ is a full subcategory of $\mathcal{B}$ such that
$G$ preserves exact sequences in $\mathcal{D}$, then for any $N\in
\mathcal{D}$, $$\elgxd(G(N))\leq\elxd(N).$$\item If
$GF(\mathcal{Y})\subseteq\mathcal{Y}$, and
$\mathcal{J}\supseteq\mathcal{Y}$ is a full subcategory of
$\mathcal{A}$ such that $F$ preserves exact sequences in
$\mathcal{J}$, then for any $M\in \mathcal{J}$,
$$\erfyd(F(M))\leq\eryd(M).$$
\end{enumerate}
\end{thm}
\proof 1. Let $\elxd(N)=m<\infty$ with $N\in\mathcal{D}$, then there
is an exact left $\X$-resolution of $N$ of length $m$ in
$\mathcal{D}$
$$0\to X_m\to\cdots\to X_1\to X_0\to N\to 0.$$
Since $G$ preserves exact sequences in $\mathcal{D}$ and $G$
preserves $\X$-precovers by \cite[Proposition 4]{GT}, we obtain an
exact left $G(\X)$-resolution of $G(N)$ $$0\to G(X_m)\to\cdots\to
G(X_1)\to G(X_0)\to G(N)\to 0.$$ So $\elgxd(G(N))\leq m$ as
required.

2. We may assume that $\eryd(M)=n<\infty$ with $M\in \mathcal{J}$.
Then there is an exact right $\Y$-resolution of $M$ of length $n$ in
$\mathcal{J}$
$$0\to M\to Y^{0}\to Y^{1}\to\cdots\to Y^{n}\to 0.$$
Since  $F$ preserves exact sequences in $\mathcal{J}$ and $F$
preserves $\Y$-preenvelopes by Lemma \ref{lem-compa-1}, we have an
exact right $F(\Y)$-resolution of $F(M)$ $$0\to F(M)\to F(Y^{0})\to
F(Y^{1})\to\cdots\to F(Y^{n})\to 0.$$ Hence $\erfyd(F(M))\leq n$ as
desired.\cqfd

We have also the following situation.

\begin{lem}\label{lem-compa2} Let $(F,G):  \mathcal{A}\rightarrow \mathcal{B}$ be an  adjoint pair with $\mathcal{X}\subseteq\mathcal{A}$ and
$\mathcal{Y}\subseteq\mathcal{B}$.
\begin{enumerate}
\item If $f: U\rightarrow V$ is an $\mathcal{X}$-precover in $\mathcal{J}$, then   $F(f): F(U)\rightarrow F(V)$ is an $F(\mathcal{X})$-precover.
\item If $\alpha: K\rightarrow H$ is a $\mathcal{Y}$-preenvelope in $\mathcal{D}$, then  $G(\alpha): G(K)\rightarrow G(H)$ is a $G(\mathcal{Y})$-preenvelope.
\end{enumerate}
\end{lem}
 \proof  1. Let $g: F(W)\rightarrow F(V)$
be a morphism with $W\in \mathcal{X}$, then there exists $h: W\rightarrow U$ such that $fh=\eta_{V}^{-1}G(g)\eta_{W}$. Thus we have
$$F(f)F(h)=\varepsilon_{F(V)}F(\eta_{V})F(fh)=\varepsilon_{F(V)}F(\eta_{V})F(\eta_{V}^{-1}G(g)\eta_{W})$$
$$=\varepsilon_{F(V)}F(\eta_{V})F(\eta_{V}^{-1})FG(g)F(\eta_{W})=\varepsilon_{F(V)}FG(g)F(\eta_{W})=g\varepsilon_{F(W)}F(\eta_{W})=g.$$
So $F(f): F(U)\rightarrow F(V)$ is  an $F(\mathcal{X})$-precover of $F(V)$.\\

2. Let $\beta: G(K)\rightarrow G(L)$ be a morphism with $L\in \mathcal{Y}$, then we get the following commutative diagram:
$$\xymatrix{G(K)\ar[d]^{\beta} \ar[r]^{\eta_{G(K)}}&GFG(K)\ar[d]^{GF(\beta)}\\
G(L)\ar[r]^{\eta_{G(L)}}&GFG(L).}$$ There exists $\gamma: H\rightarrow L$ such that $\gamma\alpha=\varepsilon_{L}F(\beta)\varepsilon_{K}^{-1}$. Therefore we have
$$G(\gamma)G(\alpha)=G(\gamma\alpha)G(\varepsilon_{K})\eta_{G(K)}=G(\varepsilon_{L}F(\beta)\varepsilon_{K}^{-1})G(\varepsilon_{K})\eta_{G(K)}$$
$$=G(\varepsilon_{L})GF(\beta)G(\varepsilon_{K}^{-1})
G(\varepsilon_{K})\eta_{G(K)}=G(\varepsilon_{L})GF(\beta)\eta_{G(K)} =G(\varepsilon_{L})\eta_{G(L)}\beta=\beta.$$ Thus $G(\alpha): G(K)\rightarrow G(H)$ is a
$G(\mathcal{Y})$-preenvelope of $G(K)$.
\cqfd

\begin{thm}\label{thm-compa2} Let $(F,G):  \mathcal{A}\rightarrow \mathcal{B}$ be an  adjoint pair with $\mathcal{X}\subseteq\mathcal{A}$ and
$\mathcal{Y}\subseteq\mathcal{B}$.
\begin{enumerate}
\item  If $\mathcal{J}\supseteq\mathcal{X}$ is a full subcategory of $\mathcal{A}$ closed under kernels of epimorphisms such that $\eta_{M}: M\rightarrow GF(M)$ is an
isomorphism for any $M\in\mathcal{J}$ and $F$ preserves exact
sequences in $\mathcal{J}$, then
$$\elfxd(F(M))\leq\elxd(M).$$
\item  If $\mathcal{D}\supseteq\mathcal{Y}$ is a full subcategory of $\mathcal{B}$ closed under cokernels of monomorphisms such that $\varepsilon_{N}:
FG(N)\rightarrow N$ is an isomorphism for any $N\in\mathcal{D}$ and
$G$ preserves exact sequences in $\mathcal{D}$, then
$$\ergyd(G(N))\leq\eryd(N).$$
\end{enumerate}
\end{thm}
\proof 1. We may assume that  $\elxd(M)=n<\infty$ with $M\in\mathcal{J}$, then there is an exact left $\X$-resolution of $M$ of length $n$ in $\mathcal{J}$
$$0\to X_n\to\cdots\to X_1\to X_0\to M\to 0.$$
Since $\mathcal{J}$ is closed under kernels of epimorphisms, $F$
preserves exact sequences and $\X$-precovers in $\mathcal{J}$, we
get an exact left $F(\X)$-resolution of $F(M)$
$$0\to F(X_n)\to\cdots\to F(X_1)\to F(X_0)\to F(M)\to 0.$$ Hence
$\elfxd(F(M))\leq n$.

2. We assume that $\eryd(N)=m<\infty$ with $N\in\mathcal{D}$. Then there is an exact right $\Y$-resolution of $N$ of length $m$ in $\mathcal{D}$
$$0\to N\to Y^{0}\to Y^{1}\to\cdots\to Y^{m}\to 0.$$
Since $\mathcal{D}$ is closed under cokernels of monomorphisms, $G$
preserves exact sequences and $\Y$-preenvelopes in $\mathcal{D}$, we
obtain an exact right $G(\Y)$-resolution of $G(N)$ $$0\to G(N)\to
G(Y^{0})\to G(Y^{1})\to\cdots\to G(Y^{m})\to 0.$$ So
$\ergyd(G(N))\leq m$. \cqfd

As a consequence we find again \cite[Theorem 2.11]{TW}.

For the following result we set $\mathcal{X}=Add_R(C)$ and
$\mathcal{Y}=Prod_R(C)$ whose objects are called, respectively,
${\mathscr P}_C$-projective  and ${\mathscr I}_C$-injective modules.
When $R$ is commutative and $C$ is a semidualizing module, Ding and
Geng in \cite{DG} showed that $Add_R(C)=\{C\otimes_{R} P: P$ is  a
projective $R$-module\} and $Prod_R(C)=\{\Hom_{R}(C,I): I$ is  an
injective $R$-module\}. So as in \cite{TW}, we use the following
notations ${\mathscr P}_C-pd_{R}(M)=\elxd(M)$ and ${\mathscr I}_C-
id_{R}(M) = \eryd(M) $ called ${\mathscr P}_C$-projective and
${\mathscr I}_C$-injective dimensions, respectively. For $C=R$ we
find $\pd_{R}(M)$  and $\id_{R}(M)$, the classical projective and
injective dimensions respectively.

\begin{cor}[\cite{TW}, Theorem 2.11]\label{cor: 4.3} Let $C$ be a semidualizing $R$-module over a commutative ring $R$. The following equalities hold.
\begin{enumerate}
\item $\pd_{R}(M) = {\mathscr P}_C-pd_{R}(C\otimes_{R} M)$.\item ${\mathscr I}_C- id_{R}(M) = \id_{R}(C\otimes_{R} M)$. \item
${\mathscr P}_C-pd_{R}(M) = \pd_{R}(\Hom_{R}(C,M))$.\item $\id_{R}(M) ={\mathscr I}_C- id_{R}(\Hom_{R}(C,M))$.
\end{enumerate}
\end{cor}
\proof Let $F=C\otimes_{R}-$ and $G=\Hom_{R}(C,-)$.

1. If ${\mathscr P}_C-pd_{R}(C\otimes_{R} M)<\infty$, then
$M\cong\Hom_{R}(C,C\otimes_{R} M)$ by \cite[Theorem 2.8 and
Corollary 2.9]{TW}. Let $\mathcal{X}={\mathscr P}_C$,
$\mathcal{D}=\{D: {\mathscr P}_C-pd_{R}(D)<\infty\}$. Then
$\pd_{R}(M)=\pd(\Hom_{R}(C,C\otimes_{R} M)) \leq
 {\mathscr
P}_C-pd_{R}(C\otimes_{R} M)$ by  Theorem \ref{thm-compa-1}(2).

Conversely, let $\pd_{R}(M)<\infty$, then
$M\cong\Hom_{R}(C,C\otimes_{R} M)$ by \cite[1.9(b)]{TW}. Let
$\mathcal{X}=Proj(R)$, $\mathcal{J}=\{J: \pd(J)<\infty\}$. Then
${\mathscr P}_C-pd_{R}(C\otimes_{R} M)\leq \pd_{R}(M)$ by Theorem
\ref{thm-compa2}(1).

It follows that $\pd_{R}(M) = {\mathscr P}_C-pd_{R}(C\otimes_{R}
M)$.

2. Suppose that ${\mathscr I}_C- id_{R}(M)<\infty$. Let
$\mathcal{Y}={\mathscr I}_C$, $\mathcal{J}=\{J: {\mathscr I}_C-
id_{R}(J)<\infty\}$. Then $\id_{R}(C\otimes_{R} M) \leq {\mathscr
I}_C- id_{R}(M)$ by Theorem  \ref{thm-compa-1}(1).

Conversely, let $\id_{R}(C\otimes_{R} M)<\infty$, then
$M\cong\Hom_{R}(C,C\otimes_{R} M)$ by \cite[1.9(b)]{TW}. Let
$\mathcal{Y}=Inj (R)$, $\mathcal{D}=\{D: \id_{R}(D)<\infty\}$. Then
${\mathscr I}_C- id_{R}(M)={\mathscr I}_C-
id_{R}(\Hom_{R}(C,C\otimes_{R} M))\leq \id_{R}(C\otimes_{R} M)$ by
Theorem \ref{thm-compa2}(2).

It follows that ${\mathscr I}_C- id_{R}(M) = \id_{R}(C\otimes_{R}
M)$.

3. Suppose that ${\mathscr P}_C-pd_{R}(M)<\infty$. Let
$\mathcal{X}={\mathscr P}_{C}$, $\mathcal{D}=\{D: {\mathscr
P}_C-pd_{R}(D)<\infty\}$. Then $\pd_{R}(\Hom_{R}(C,M))\leq{\mathscr
P}_C-pd_{R}(M)$ by Theorem \ref{thm-compa-1}(2).

Conversely, let $\pd_{R}(\Hom_{R}(C,M))<\infty$, then $M\cong
C\otimes_{R}\Hom_{R}(C, M)$. Let $\mathcal{X}=Proj(R)$,
$\mathcal{J}=\{J: \pd_{R}(J)<\infty\}$. Then  ${\mathscr
P}_C-pd_{R}(M)={\mathscr P}_C-pd_{R}(C\otimes_{R}\Hom_{R}(C, M))\leq
\pd_{R}(\Hom_{R}(C,M))$ by Theorem \ref{thm-compa2}(1).

It follows that ${\mathscr P}_C-pd_{R}(M) = \pd_{R}(\Hom_{R}(C,M))$.

4. If ${\mathscr I}_C- id_{R}(\Hom_{R}(C,M))<\infty$, then $M\cong
C\otimes_{R}\Hom_{R}(C, M)$. Let $\mathcal{Y}={\mathscr I}_C$,
$\mathcal{J}=\{J: {\mathscr I}_C- id_{R}(J)<\infty\}$. Then
$\id_{R}(M)=\id( C\otimes_{R}\Hom_{R}(C, M)) \leq{\mathscr I}_C-
id_{R}(\Hom_{R}(C,M))$ by  Theorem  \ref{thm-compa-1}(1).

Conversely, let $\id_{R}(M)<\infty$, then $M\cong
C\otimes_{R}\Hom_{R}(C, M)$.  Let $\mathcal{Y}=Inj (R)$,
$\mathcal{D}=\{D: \id(D)<\infty\}$. Then ${\mathscr I}_C-
id_{R}(\Hom_{R}(C,M))\leq \id_{R}(M)$ by Theorem
\ref{thm-compa2}(2).

It follows that $\id_{R}(M) ={\mathscr I}_C- id_{R}(\Hom_{R}(C,M))$.
\cqfd

\section{Global homological dimensions}
\hskip .5cm  In this section, $\X$ and $\Y$  will be, respectively, a   precovering  and a preenveloping  subcategories of $\A$ which are self-orthogonal
and closed under extensions  and direct summands.

The aim of this section is to study the  global dimensions relative
to a subcategory.

Let $\ZZ$ be  a subcategory of  $\A$ which has enough projectives.  The left global projective dimension of   $\ZZ$, $\lgldim(\ZZ)$,  is the supremum of
 the set of projective dimensions of all objects of $\ZZ$. When $\ZZ$ contains $\X$, we define the left $\X$-global projective dimension of
 $\ZZ$, $\lXgldim(\ZZ)$,  as the supremum of the set of  $\elxd(M)$   of all objects $M$  of $\ZZ$.

We first consider the global case $\ZZ=\A$.

\begin{prop} \label{prop-gldim} Assume  $\A$ to have enough projectives.  Then, for a   positive integer $n$,  $\lXgldim(\A)= n$ if and only if  $\lgldim(\A)= n$ and
$\X$ coincides with $P(\A)$ the  subcategory of all projective
objects of $\A$.
\end{prop}
\proof $ \Rightarrow $. Assume that $\lXgldim(\A)\leq n$. Then, clearly $\X$
contains $P(\A)$. Now consider an object $C$ of $\X$ and   a short
exact sequence in $\A$, $0\rightarrow K \rightarrow P \rightarrow C
\rightarrow 0$, where $P$ is projective. By hypothesis and
Proposition \ref{prop-elxd}, $\Ext^{\geq 1}(C,K) = 0$. This implies
that the
sequence splits and so $C$ is projective.\\
$ \Leftarrow $. Obvious.\cqfd

As an application we get the following generalization of
\cite[Proposition 5.1]{TW} using different arguments.

Let us denote l.$Add(C)\! - \!$ gl.dim($\ZZ$) by $\lPCgldim(\ZZ)$.
In particular, when $\ZZ=R$-$\Mod$,   $\lPCgldim(\ZZ)$ will be
simply denoted by $\lPCgldim(R)$.

An $R$-module $C$ is  said to be $\Sigma$-self-orthogonal if
$\Ext^{\geq 1}_R (C,C^{(I)}) = 0$ for every set $I$ (that is,
$C^{(I)}\in Add_R(C)^{\perp}$ for every set $I$).

\begin{cor}  \label{cor-gldim-2}
Let $C$ be a    $\sum$-self-orthogonal $R$-module.  Then, for every
positive integer $n$,  $\lPCgldim(R)= n$ if and only if $\lgldim(R)=
n$ and $C$ is  a projective generator of $R$-$\Mod$, that is,
$Add_R(C)=P(R)$.

In particular,  $\lPCgldim(R)= 0$ if and only if $R$-$\Mod =
Add_R(C)=P(R)$.
\end{cor}

Note that the condition on $C$ to be a generator of $R$-Mod  cannot
be dropped. Indeed, take a semisimple ring $R=k_1\times \cdots
\times k_n$, where $k_i$ is  a field for each $i$. We have
$\lgldim(R)=0$ (i.e.,  $R$-Mod=$P(R)$)  and, for example, $C=k_1$ is
a $\sum$-self-orthogonal $R$-module. But $Add_R(C)\neq P(R)$.

If we consider the  ${\mathscr P}_C$-projective dimension of only finitely generated modules we get the following result.

\begin{prop}  \label{pro-gldim-2-fg}
Let $C$ be a finitely presented and  $\sum$-self-orthogonal $R$-module.  Then, for a   positive integer $n$,   the  ${\mathscr P}_C$-projective dimension of
every finitely generated module is at
most $  n$ if and only if $\lgldim(R)\leq  n $ and $C$ is  a projective generator of $R$-\Mod.
\end{prop}
\proof
 $ \Rightarrow $. By hypothesis we can show that $R$ is ${\mathscr P}_C$-projective, and so    $Add_R(C)$ contains
$P(R)$. On the other hand,  $C$ is projective. Indeed, since $C$ is finitely presented, there is    a short exact sequence
 of $R$-module $0\rightarrow K \rightarrow P \rightarrow C \rightarrow 0$,
where $P$ is a finitely generated projective  $R$-module and $K$ is finitely generated. By hypothesis and Proposition \ref{prop-elxd}, $\Ext^{\geq 1}_R (C,K) = 0$.
This implies that the
sequence splits and so $C$ is projective. Therefore,  $C$ is  a projective generator of $R$-\Mod, as desired.  \\
$ \Leftarrow $. Obvious.\cqfd

  Notice that even if we do not   assume that $C$ is finitely presented in Proposition \ref{pro-gldim-2-fg} above,   the condition "the  ${\mathscr P}_C$-projective dimension
of  every finitely generated module is at most a  positive integer
$n$" implies that the injective dimension of  every  ${\mathscr
P}_C$-projective module is at most $n$. Then,  if we further assume
that   $R$ is Noetherian, then $R$ will be $n$-Gorenstein  (since
$R$ is ${\mathscr P}_C$-projective). Then, every finitely generated
module has  projective dimension at most $  n$, and so
$\lgldim(R)\leq  n $. This leads to the following question: whether
the condition ``$C$ is finitely presented" in Proposition
\ref{pro-gldim-2-fg} can be removed?

If we discuss the case $n=0$, we can show that we have a positive
answer to this question. Indeed, being a finitely generated
$R$-module, $R$  is ${\mathscr P}_C$-projective and so every
projective $R$-module is ${\mathscr P}_C$-projective. Then, by the
argument above,  every  projective $R$-module is injective. This
implies that $R$ is quasi-Frobenius, in particular, Noetherian.
Then, as above we conclude that $R$ is semisimple.

Similarly we define, when    $\A$ has enough injectives, the right
global injective dimension of   $\A$, $\rgldim(\A)$,  as the
supremum of the set of   injective dimensions of all objects. The
right $\Y$-global injective dimension of   $\A$, $\rYgldim(\A)$, is
the supremum of the set of  $\eryd(M)$   of all objects $M$.

As a dual result of Proposition \ref{prop-gldim} we have the following result.

\begin{prop} \label{cor-gldim-dual} Assume  $\A$ to have enough injectives.  Then, for a   positive integer $n$,  $\rYgldim(\A)= n$ if and only if  $\rgldim(\A)= n$
and $\Y$ coincides with $I(\A)$
the  subcategory of all injective objects of $\A$.
\end{prop}

As a consequence, we get the following result.

 Let us denote r.$Prod(C)\! - \!$ gl.dim($\ZZ$) by $\rICgldim(\ZZ)$. In particular,
when $\ZZ=R$-$\Mod$,   $\rICgldim(\ZZ)$ will be simply denoted by
$\rICgldim(R)$.

An $R$-module $C$ is  said to be $\prod$-self-orthogonal if
$\Ext^{\geq 1}_R (C,C^{I}) = 0$ for every set $I$.

\begin{cor}  \label{cor-gldim-2-dual}
Let $C$ be a   $\prod$-self-orthogonal $R$-module.  Then, for a
positive integer $n$,  $\rICgldim(R)= n$ if and only if $\rgldim(R)=
n$ and $C$ is  an injective cogenerator of $R$-$\Mod$, that is,
$Prod_R(C)=I(R)$.

In particular,  $\rICgldim(R)= 0$ if and only if   $R$-$\Mod= Prod_R(C)=I(R)$.
\end{cor}

In \cite[Corollary 4.4]{BGO2}, it was proved that, if $C$ is a
$\sum$-self-orthogonal, self-small and Hom-faithful $R$-module, then
$\lgldim (S)=\lPCgldim (R)$, where   $S=\End_R(C)$, the endomorphism
ring of $C$.  This result with Corollary \ref{cor-gldim-2} show
that, under the condition above on $C$, if $\lPCgldim(R)$ is finite
then,   $\lgldim (S)=\lPCgldim (R)$. There is a classical result by
Miyashita \cite{M86}  that relates $\lgldim (S)$, $\lPCgldim (R)$
and the projective dimension of $C$ when $C$ is assumed to be a
tilting module (see \cite[Corollary  of  Proposition  2.4]{M86}).
 In \cite[Theorem 1.3]{T97},  Trlifaj related the  $\lgldim (S)$ and $\lPCgldim (R)$ when $C$ is a $*$-module (see also \cite{W5,W6}).  Then,
 a natural question arises: How are   $\lgldim (S)$ and $\lPCgldim (R)$
 related when $C$ is a $\sum$-self-orthogonal, self-small, and Hom-faithful $R$-module?

One can be interested in studying l.$Add(C)\! - \!$gl.dim$(\ZZ)$
for some interesting particular cases of subcategories. For
instance, $\sigma [C]$ and $\textrm{Gen}[C]$. Recall that an
$R$-module  $N$ is said to be $C$-generated if there exists an exact
sequence of $R$-modules
$$0 \to K \to C^{(\Lambda)} \to N \to 0\ \textrm{for\  some\  set\ } \Lambda.$$

$\textrm{Gen}[C]$ (resp., $\sigma [C]$) will denote the full
subcategories of $R$-Mod whose objects are $C$-generated (resp.,
submodules of $C$-generated modules) \cite{Wis91}.

\begin{prop} \label{prop-gldim-sigma}
Let $C$ be a   $\Sigma$-self-orthogonal $R$-module. The following
assertions hold:
\begin{enumerate}
   \item  $l.Add(C)\!-\!gl.dim(\sigma [C])=0$  if and only if  $C$ is a semisimple module.  \item $l.Add(C)\!-\!gl.dim(\sigma [C])\leq 1$ if and only if $Add(C)$ is closed under submodules.
\end{enumerate}
\end{prop}
\proof 1. Easy to prove.\\
2. $\Rightarrow.$ Let $M$ be a submodule of a module $N\in Add(C)$.
Then,   $M\in \sigma [C]$. Since l.$Add(C)\! - \!$gl.dim$(\sigma
[C])\leq 1$, $M\in Add(C)^{\perp}$ (by Proposition \ref{prop-elxd}).
This shows that the sequence  $0\to M\to N\to M/N\to  0$ is
$\Hom_R(Add(C),-)$-exact. On the other hand,     $N/M\in \sigma
[C]$, then ${\mathscr P}_C-pd_{R}(N/M)\leq 1$.
So, by Proposition \ref{prop-elxd}, $M\in Add(C)$, as desired.\\
$\Leftarrow.$ Consider an $R$-module $M\in \sigma [C]$. Then,    $M$
is a submodule of a module $N \in  \textrm{Gen}[C]$.  There exists a
short exact sequence of $R$-modules of the following form
$$0 \to K \to C^{(\Lambda)} \to N \to 0 \ \textrm{for\  some\  set\ } \Lambda.$$
Then we get the following pullback diagram
$$   \xymatrix{
     &  0 \ar[d]  & 0 \ar[d]  &   &  \\
 &  K\ar[d] \ar@{=}[r]&  K \ar[d]  &   &  \\
 0\ar[r]& D\ar@{-->}[d] \ar@{-->}[r] & C^{(\Lambda)} \ar[d] \ar[r] &  N/M \ar@{=}[d]  \ar[r] & 0\\
0\ar[r]& M \ar[r] \ar[d]&  N  \ar[d] \ar[r]& N/M\ar[r] & 0\\
 & 0 &0  &   & }  $$
By hypothesis $D$ and $K$ are in $Add(C)$, as desired.\cqfd

Now we give a result on  l.$Add(C)\! - \!$gl.dim(\textrm{Gen}$[C])$.

We will say that a submodule $N$ of an $R$-module $M$ is $C$-pure if the sequence $0\rightarrow N\rightarrow M\rightarrow M/N\rightarrow 0$ is $\Hom_R(C,-)$-exact.

\begin{prop} \label{prop-gldim-Gen}
Let $C$ be a   $\Sigma$-self-orthogonal $R$-module.
The following assertions hold:
\begin{enumerate}
   \item  $l.Add(C)\!-\!gl.dim(\textrm{Gen} [C])=0$  if and only if, for every module $N\in Add(C)$, every $C$-pure submodule of $N$ is a direct summand.
\item $l.Add(C)\!-\!gl.dim(\textrm{Gen} [C])\leq 1$ if and only if $Add(C)$ is closed under $C$-pure submodules.
\end{enumerate}
\end{prop}
\proof 1. $\Rightarrow.$ Let $M$ be a $C$-pure submodule of a module
$N\in Add(C)$. Then, $N/M\in$ Gen$[C]=Add(C)$. Hence the sequence
$0\rightarrow M\rightarrow N\rightarrow N/M\rightarrow 0$ splits
and so $M$ is a direct summand of $N$.\\
$\Leftarrow.$ Let $N\in$ Gen$[C]$ and $0\rightarrow K\rightarrow
G\rightarrow N\rightarrow 0$ be an exact sequence with $G\rightarrow
N$ an epic $Add(C)$-precover. Then $K$ is a $C$-pure submodule of
$G$, hence by hypothesis, the sequence splits and so $N\in Add(C)$.

2. $\Rightarrow.$ Let $M$ be a $C$-pure submodule of a module $N\in
Add(C)$. Then, $N/M\in$ Gen$[C]$. We consider the exact sequence
$0\rightarrow C^1\rightarrow C^0\rightarrow N/M\rightarrow 0$ with
$C^1, C^0\in Add(C)$.  On the other hand, using the fact that  $C$
is  a $\Sigma$-self-orthogonal $R$-module and that $M$ is a $C$-pure
submodule of a module $N\in Add(C)$, we deduce that $\Ext^1(D,M)=0$
for all $D   \in Add(C)$. By the relative version of
Schanuel's lemma \cite[Lemma 8.6.3]{EJb}, $C^1\oplus M\cong C^0\oplus N$, hence $M\in Add(C)$.\\
$\Leftarrow.$ Consider an $R$-module $M\in$ Gen$[C]$. Then, there is
an exact sequence $0\rightarrow K\rightarrow D\rightarrow
M\rightarrow 0$ where $D\rightarrow M$ is an epic $Add(C)$-precover.
So $K\in Add(C)$ by hypothesis and hence we have the result. \cqfd


\section{Relative (co)homology groups}
\hskip .5cm In this  section, $\X$ and $\X'$ (resp., $\Y$ and $\Y'$), when
considered  as  subcategories of   abelian categories, will be
precovering (resp., preenveloping), self-orthogonal  and    closed
under extensions  and direct summands. We will be interested in
relative (co)homology groups. The results established in this section are inspired from the ones of \cite{SSW} and \cite{TW}.\\

Recall that the relative Ext-groups with respect to   $\mathcal{X}$ are defined as
$$\Ext^i_\mathcal{X}(M,N)=H^i\Hom(X,N),$$
where $X$ is a left proper $\mathcal{X}$-resolution of $M$ and $i\geq 0$.
Analogously, for the preenveloping class $\mathcal{Y}$, we have
$$\Ext^i_\mathcal{Y}(M,N)=H^i\Hom(M,Y),$$
where $Y$ is a right coproper $\mathcal{Y}$-resolution of $N$ and
$i\geq 0$.

The first main result of this section generalizes \cite[Theorem
3.2]{TW}. It characterizes relative dimensions using relative
cohomology groups.

 It needs the   following lemma.

\begin{lem}\label{lem-rel-func}
 If $\phi: X\rightarrow M$ is an $\X$-precover of an object $M$, then $\phi': X\rightarrow M'=Im(\phi)$ defined by $\phi'(x)= \phi(x)$ for all $x\in X$ is an epic $\X$-precover of $M'$ and
 $\Ext^{1}(N,\Ker(\phi))=0$ for every object $N\in \X$.
\end{lem}
\proof Since $\X$ is self-orthogonal, we only need to prove that
$\phi'$ is an epic $\X$-precover of $M'$. For this, consider a
homomorphism $f:X'\to M'$, where $X'$ is an object of $\X$. Then,
with the injection $i:M' \to M$, we get the following homomorphism
$if:X'\to M$. Since $\phi$ is an  $\X$-precover, there exists $f':
X'\to X$ such that $\phi f'=if$. Since $i \phi'=\phi$, we get $i
\phi' f'=if$. This implies that $  \phi' f'= f$ since $i$ is monic,
as desired.\cqfd

\begin{thm}\label{thm-rel-func}
Assume $\X$ to be  $\Hom$-faithful, then for an integer $ n\geq 0$ the following are equivalent.
\begin{enumerate}
\item $\lxd( M)\leq n$.
\item $\Ext^{n+1}_{\X}(M,N)=0$ for every object $N\in \A$.
\item $\Ext^{i}_{\X}(M,N)=0$ for every $i\geq n+1$ and every object $N\in \A$.
\end{enumerate}
\end{thm}

\proof
$1.\Rightarrow 3.$ and $3.\Rightarrow 2.$ are clear. \\
$2.\Rightarrow 1$. An inductive argument shows that we  only need to
prove the result for   the case $n=0$. Then assume that
$\Ext^{1}_{\X}(M,N)=0$ for every object $N\in \A$. Consider a left
proper $\X$-resolution of $M$
$$\xymatrix @!0 @R=10mm  @C=1.5cm  { \cdots \ar[r]&  X_2 \ar[dr]_{f_2} \ar[rr]^{d_2=i_1f_2} &&  X_1  \ar[dr]_{f_1} \ar[rr]^{d_1=i_0f_1} &&  X_0   \ar[r]^{d_0} & M \ar[r] & 0  \\
  & & K_1 \ar[ur]_{i_1} &  &   K_0 \ar[ur]_{i_0}  &       \\
   & 0 \ar[ur]  & &   0 \ar[ur]   &       &  & }$$
 We have in particular  $\Ext^{1}_{\X}(M,K_0)=0$.
Then $\ker(\Hom(d_2,K_0))=\Im(\Hom(d_1,K_0))$. Since
 $$\ker(\Hom(d_2,K_0))= \{ g\in \Hom(X_1,K_0) \ / \ gi_1f_2=0  \}$$ and
$$ \Im(\Hom(d_1,K_0))= \{ hi_0f_1 \in \Hom(X_1,K_0)\ / \ h\in \Hom(X_0,K_0) \},$$
we deduce that, for every $g \in \Hom(X_1,K_0)$ such that
$gi_1f_2=0$ there exists $h\in \Hom(X_0,K_0)$ such that $g=hi_0f_1$.
In particular,  there exists $h\in \Hom(X_0,K_0)$ such that
$f_1=hi_0f_1$ (since   $f_1i_1f_2=0$). By applying $\Hom(X,-)$ and
setting for $X\in \X$,  $f^*=\Hom(X,f)$ we get
  $ f^*_1=h^*i^*_0f^*_1 $. Then $h^*i^*_0= id_{\Hom(X,K_0)}$ since $f^*_1$ is epic. Now consider the exact sequence:
$$
0  \rightarrow K \rightarrow   K_{0}  \stackrel{hi_0}{\rightarrow} K_0  \rightarrow \coker(hi_0) \rightarrow  0.
$$
We have then   $ \Hom(X,K)=0 $ for every $X\in \X$ and this implies
that $K=0$ and hence $hi_0$ is monic. Now, by Lemma
\ref{lem-rel-func},  $\Ext^{1}(X,K_0 )=0$ for every $X\in \X$. Then,
we get the following   exact sequence:
$$ 0\rightarrow \Hom(X,K_0)   \rightarrow \Hom(X,K_0)   \rightarrow  \Hom(X, \coker(hi_0)) \rightarrow 0.
$$
Thus $ \Hom(X, \coker(hi_0)) =0$ because $h^*i^*_0= id_{K_0}$.
Therefore $hi_0$ is an automorphism, which means that  the exact
sequence $0\rightarrow K_0 \rightarrow  X_0 \rightarrow M'  =\Im(
X_0 \rightarrow M  )\rightarrow 0$ splits and then $M'\in \X$. It
follows that the homomorphism $M'=X_0/K_0 \rightarrow M$ induced by
$X_0 \rightarrow M$ is a monic $\X$-precover of $M$, and thus  it is
an isomorphism by Theorem  \ref{thm-eld-ld}, as desired.\cqfd

It is worth noting that in \cite[Section 3]{H} Holm studied the
relation between the relative cohomology group and the relative
dimension. For that reason he introduced the notion of  almost
epimorphisms which serves as a sufficient condition to establish a
result analogue to Theorem \ref{thm-rel-func}.

Our second aim is the   balance results for relative (co)homology groups. First  we give the following   result      which generalizes   \cite[Theorem 4.1]{TW}.

\begin{thm}\label{thm-compa-3} Let $(F,G):  \mathcal{A}\rightarrow \mathcal{B}$ be an
adjoint pair with $\mathcal{Y}\subseteq\mathcal{A}$ and
$\mathcal{X}\subseteq\mathcal{B}$.
\begin{enumerate}\item If $FG(X)\cong X$ for any $X\in\mathcal{X}$, then $\Ext_{\mathcal{X}}^{i}(M,N)\cong\Ext_{G(\mathcal{X})}^{i}(G(M),G(N))$ for
any $M,N\in\mathcal{B}$.
\item  If $GF(Y)\cong Y$ for any $Y\in\mathcal{Y}$, then $\Ext_{\mathcal{Y}}^{i}(M,N)\cong\Ext_{F(\mathcal{Y})}^{i}(F(M),F(N))$ for any $M,N\in\mathcal{A}$.
\end{enumerate}
\end{thm}
\proof 1. Let $\mathcal{Z}^{+}: \cdots\to X_1\to X_0\to M\to 0$ be a
left proper $\X$-resolution of $M$ and $\mathcal{Z}: \cdots\to
X_1\to X_0\to 0$ be the deleted complex. Then $G(\mathcal{Z}^{+}):
\cdots\to G(X_1)\to G(X_0)\to G(M)\to 0$ is a left proper
$G(\X)$-resolution of $G(M)$ by \cite[Proposition 4]{GT}. So
$\Ext_{G(\mathcal{X})}^{i}(G(M),G(N))=H^{i}\Hom(G(\mathcal{Z}),G(N))\cong
H^{i}\Hom(FG(\mathcal{Z}),N)\cong
H^{i}\Hom(\mathcal{Z},N)=\Ext_{\mathcal{X}}^{i}(M,N).$

2. Let $\mathcal{W}^{+}:  0\to N\to Y^{0}\to Y^{1}\to\cdots$ be a
right coproper  $\Y$-resolution of $N$ and $\mathcal{W}: 0\to
Y^{0}\to Y^{1}\to\cdots$ be the deleted complex. Then
$F(\mathcal{W}^{+}): 0\to F(N)\to F(Y^{0})\to F(Y^{1})\to\cdots$ is
a right coproper $F(\Y)$-resolution of $F(N)$ by the proof of
Theorem \ref{thm-compa-1}(1). So we have
$\Ext_{F(\mathcal{Y})}^{i}(F(M),F(N))=H^{i}\Hom(F(M),F(\mathcal{W}))\cong
H^{i}\Hom(M,GF(\mathcal{W}))\cong
H^{i}\Hom(M,\mathcal{W})=\Ext_{\mathcal{Y}}^{i}(M,N).$ \cqfd

\medskip

As an immediate consequence of Theorem \ref{thm-compa-3}, we get the following  result.

\begin{cor}[\cite{TW}, Theorem 4.1] \label{cor: 4.5} Let $C$ be a semidualizing $R$-module over a commutative ring $R$ and let $M$ and $N$ be $R$-modules.
There exist isomorphisms:
\begin{enumerate}\item
$\Ext_{{\mathscr
P}_C}^{i}(M,N)\cong\Ext_{R}^{i}(\Hom_{R}(C,M),\Hom_{R}(C,N))$.\item
$\Ext_{{\mathscr
I}_C}^{i}(M,N)\cong\Ext_{R}^{i}(C\otimes_{R}M,C\otimes_{R}N)$.
\end{enumerate}
\end{cor}

Let $\mathcal{E},\; \mathcal{F} \subseteq {\cal A}$ be two classes
of objects. We say that the couple $(\mathcal{E}, \mathcal{F})$ is a
balanced pair if $\mathcal{E}$ is precovering, $\mathcal{F}$ is
preenveloping and the bifunctor $\Hom_{\cal A}(-,-)$  is right
balanced by $\mathcal{E}\times \mathcal{F}$ in the sense of
\cite[Definition 8.2.13]{EJb}.

Note that in case where $(\mathcal{E}, \mathcal{F})$ is a balanced
pair we have, for two objects $M$ and $N$, the natural isomorphisms
(see \cite[Theorem 8.2.14]{EJb}):
$$\Ext^i_\mathcal{E}(M,N)\cong \Ext^i_\mathcal{F}(M,N)\cong
H^iTot(\Hom(\textbf{E},\textbf{F})),$$ where   $\textbf{E}$ is a
left proper  $\mathcal{E}$-resolution of  $M$, $\textbf{F}$ is a
right coproper  $\mathcal{F}$-resolution of  $N$ and
$Tot(\Hom(\textbf{E},\textbf{F}))$ is the total complex of the
bicomplex $\Hom(\textbf{E},\textbf{F})$.

\begin{thm}\label{thm-compa-5} Let $(F,G):  \mathcal{A}\rightarrow \mathcal{B}$ be an
adjoint pair with $\mathcal{X},\,\mathcal{Y}\subseteq\mathcal{A}$,
$\mathcal{X'},\,\mathcal{Y'}\subseteq\mathcal{B}$ subcategories. Assume $(\mathcal{X}, \mathcal{Y})$ and
$(\mathcal{X'}, \mathcal{Y'})$ are balanced pairs and $FG(X)\cong X$
naturally for all $ X\in \mathcal{X}$, $GF(Y)\cong Y$ naturally
 for all $ Y\in \mathcal{Y}$.
\begin{enumerate}\item If $G(\mathcal{X'})\subseteq
\mathcal{X}$, $FG(M)\cong M$ and $FG(N)\cong N$ naturally,  then
$$\Ext_{\mathcal{Y'}}^{i}(M,N)\cong\Ext_{F(\mathcal{Y})}^{i}(M,N).$$
\item If $F(\mathcal{Y})\subseteq \mathcal{Y'}$, $GF(L)\cong L$ and $GF(T)\cong T$ naturally, then
 $$\Ext_{\mathcal{X}}^{i}(L,T)\cong\Ext_{G(\mathcal{X'})}^{i}(L,T).$$
\end{enumerate}
\end{thm}
\proof 1. Let $X'$ be a left proper  $\mathcal{X'}$-resolution of $M$. Then
by \cite[Proposition 4]{GT}, $G(X')$ is a  left
proper  $\mathcal{X}$-resolution of $G(M)$. By hypothesis, $X'$ can be
written as $X'=F(G(X'))$ and $G(X')$ is a  left
proper $\mathcal{X}$-resolution of $G(M)$. Hence
 $\Ext_{\mathcal{Y'}}^{i}(M,N)\cong\Ext_{\mathcal{X'}}^{i}(M,N)\cong
H^i\Hom(F(G(X')), N)\cong H^iHom(G(X'),G(N))\cong
\Ext^i_{\mathcal{X}}(G(M),G(N))\cong\Ext^i_{\mathcal{Y}}(G(M),G(N))\cong^{(*)}$\linebreak$\Ext^i_{F(\mathcal{Y})}(FG(M),FG(N))\cong
\Ext^i_{F(\mathcal{Y})}(M,N)$, where the isomorphism $(*)$ holds by
Theorem \ref{thm-compa-3}(2).

2. Let $Y$ be a right coproper  $\mathcal{Y}$-resolution of $T$. Then by
\cite[Proposition 2.5]{GT1}, $F(Y)$ is  a right coproper
$\mathcal{Y'}$-resolution of $F(T)$. By hypothesis, $Y$ can be
written as $Y\cong G(F(Y))$ and $F(Y)$ is a  right coproper
$\mathcal{Y}$-resolution of $F(T)$. Hence $\Ext_{\mathcal{X}}^{i}(L,T)\cong\Ext_{\mathcal{Y}}^{i}(L,T)\cong H^i\Hom(L, G(F(Y)))\cong
H^i\Hom(F(L),F(Y))\cong  \Ext^i_{\mathcal{Y'}}(F(L),F(T))\cong
\Ext^i_{\mathcal{X'}}(F(L),F(T))\cong^{(*)}
\Ext^i_{G(\mathcal{X'})}(GF(L),GF(T))\cong $\linebreak$\Ext^i_{G(\mathcal{X'})}(L,T)$,
where the isomorphism $(*)$ holds by Theorem \ref{thm-compa-3}(1).
\cqfd

We end the paper with a result which generalizes \cite[Proposition 4.2]{TW}. It needs the following lemma.

\begin{lem}\label{thm-compa-6}Let $(F,G):  \mathcal{A}\rightarrow \mathcal{B}$ be an
adjoint pair with $\mathcal{X}\subseteq\mathcal{A}$,
$\mathcal{Y}\subseteq\mathcal{B}$ self-orthogonal subcategories.
Assume   $\mathcal{Y}\subseteq {\cal V}=\{ Y\in {\cal A}\mid
FG(Y)\cong Y\; naturally \}$ and $\mathcal{X}\subseteq {\cal W}=\{
X\in {\cal B}\mid GF(X)\cong X\; naturally\}$. The following
conditions are equivalent.
\begin{enumerate}
    \item  $(F(\mathcal{X}), \mathcal{Y})$ is balanced in ${\cal V}\times
{\cal V}$.

  \item  $(\mathcal{X}, G(\mathcal{Y}))$ is balanced in ${\cal W}\times
{\cal W}$.
\end{enumerate}
\end{lem}
\proof $1.\Rightarrow 2.$ Let $M,N\in {\cal W}$. First note that by
the equality $\varepsilon_{F(M)}\circ F(\eta_M)=id_{F(M)}$, $F({\cal
W})\subseteq {\cal V}$.

By hypothesis there is a $\Hom(-,\mathcal{Y})$-exact sequence
$$\cdots \rightarrow F(X^1)\rightarrow F(X^0)\rightarrow
F(M)\rightarrow 0.$$ Hence the sequence $$\cdots \rightarrow X^1\cong GF(X^1)\rightarrow X^0\cong GF(X^0)\rightarrow M\cong GF(M)\rightarrow 0$$
is $\Hom(-, G(\mathcal{Y}))$-exact. Analogously, let
$$0\rightarrow F(N)\rightarrow Y^0\rightarrow Y^1\rightarrow \cdots$$ be $\Hom(F(\mathcal{X}), -)$-exact. Then the sequence
$$0\rightarrow N\cong GF(N)\rightarrow G(Y^0)\rightarrow
G(Y^1)\rightarrow \cdots$$ is $\Hom(\mathcal{X}, -)$-exact.

$2.\Rightarrow 1.$ follows from the same ideas that the reverse
implication. \cqfd

\begin{prop}\label{thm-compa-7}Let $(F,G):  \mathcal{A}\rightarrow \mathcal{B}$ be an
adjoint pair with $\mathcal{X}\subseteq\mathcal{A}$,
$\mathcal{Y}\subseteq\mathcal{B}$  subcategories.
\begin{enumerate}\item
Assume that for some $M, N\in  \mathcal{B}$, $FG(M)\cong M$,
$\Ext^i_{\mathcal{X}}(G(M),G(Y))=0$, for all $ Y\in \mathcal{Y}$ and
all $ i\geq 1$ and $\Ext^i_{\mathcal{Y}}(F(X),N)=0$, for all  $ X\in
\mathcal{X}$ and all $  i\geq 1$. Then,
$\Ext^i_{F(\mathcal{X})}(M,N)\cong \Ext^i_{\mathcal{Y}}(M,N)$ for
all $ i\geq 1$.

\item Assume that for some $M, N\in  \mathcal{A}$, $GF(N)\cong N$,
$\Ext^i_{\mathcal{X}}(M,G(Y))=0$ for all $  Y\in \mathcal{Y}$ and
all $  i\geq 1$ and $\Ext^i_{\mathcal{Y}}(F(X),F(N))=0$ for all $
X\in \mathcal{X}$  and all $ i\geq 1$. Then,
$\Ext^i_{G(\mathcal{Y})}(M,N)\cong \Ext^i_{\mathcal{X}}(M,N)$ for
all $ i\geq 1$.
\end{enumerate}
\end{prop} \proof It holds easily by
\cite[Theorem 8.2.14]{EJb}. \cqfd
\bigskip

\noindent\textbf{Acknowledgment.} The second and fourth authors were
partially supported by the grant MTM2014-54439-P from Ministerio de
Econom\'{\i}a y Competitividad. The third author was partially
supported by NSFC (11771202). A part of this paper was presented in
the workshop  "Homological and homotopical tools in category theory
with applications to Algebraic Geometry, Representation Theory and
Module Theory" held in Hangzhou, China (May 2018). The first author
would like to thank the organizers for the warm hospitality and the
fruitful discussions. The authors want to express their gratitude to
Professor Nanqing Ding and the referee for the very helpful comments
and suggestions.

 Driss Bennis:   CeReMAR Center; Faculty of Sciences  B.P. 1014, Mohammed V University in Rabat, Rabat, Morocco.

\noindent e-mail address: driss.bennis@um5.ac.ma; driss$\_$bennis@hotmail.com

J. R. Garc\'{\i}a Rozas: Departamento de  Matem\'{a}ticas,
Universidad de Almer\'{i}a, 04071 Almer\'{i}a, Spain.

\noindent e-mail address: jrgrozas@ual.es

Lixin Mao: Department of Mathematics and Physics, Nanjing Institute
of Technology, Nanjing 211167, China

\noindent e-mail address: maolx2@hotmail.com

Luis Oyonarte: Departamento de  Matem\'{a}ticas, Universidad de
Almer\'{i}a, 04071 Almer\'{i}a, Spain.

\noindent e-mail address: oyonarte@ual.es

\end{document}